\documentclass[36pt, reqno]{amsart}

\usepackage{amsmath,amssymb,amscd,amsfonts,latexsym,epsfig, subfigure}

\newcommand\PP{\mathbb P} 
\newcommand\ZZ{\mathbb Z}
\newcommand\NN{\mathbb N}
\newcommand\QQ{\mathbb Q}
\newcommand{\NI}{{\noindent}}

\renewcommand\L{\mathcal L}
\renewcommand\O{\mathcal O}

\makeatletter \@addtoreset{equation}{section} \makeatother

\newtheorem{thm}[equation]{Theorem} \newtheorem{lem}[equation]{Lemma}
\newtheorem*{eg*}{Example}
\newtheorem{cor}[equation]{Corollary}
\newtheorem*{cor*}{Corollary}

\theoremstyle{definition} \newtheorem{defn}[equation]{Definition}

\title[Unstable Kodaira Fibrations]{Unstable Kodaira Fibrations}

\author{Yujen Shu}

\begin{document}
\bibliographystyle{plain}

\begin{abstract}
  Being inspired by Ross' construction of unstable products of certain
  smooth curves, we
  show that the product $C\times C$ of every smooth curve $C$ of
  genus at least 2 is not slope semistable with respect to certain polarisations.
  Besides, we
   produce examples of Kodaira-fibred surfaces of nonzero signature,
   which are not slope semistable with
  respect to some polarisations, and so they admit K\"ahler classes
  that do not contain any constant scalar curvature K\"ahler metrics.
\end{abstract}

\maketitle
\section{Introduction}
The famous Calabi-Yau \cite{ca}\cite{yau} theorem has told us that every compact complex
manifold $X$ with negative first Chern class admits a
K\"{a}hler-Einstein metric in the class $-c_1(X)$, which has a
negative constant scalar curvature. By a deformation argument due to
LeBrun and Simanca \cite{ls}, there exists a constant scalar
curvature K\"{a}hler(cscK) metric in every class near $-c_1(X)$. For
a Kodaira-fibred surface $\pi:X\to B$, Fine \cite{fine} found the
existence of cscK metrics in the classes of $-c_1(X)-rc_1(B)$, which
is far from $-c_1(X$), for large $r$ via an adiabatic limit and the
inverse theorem of Banach space. It had been an open problem that
whether every K\"{a}hler class of a compact complex manifold $X$
with negative first Chern class contains a cscK metric until
Ross \cite{ross} constructed the first example by products of curves
which fails to have any cscK metric in some K\"ahler classes.

In Ross' example \cite{ross}, he shows that the product $C\times C$ is not
semistable with respect to certain polarisations if $C$ is a simple
branched cover of $\PP_1$ of degree $k$, where $2\leq k-1<\sqrt{
\text{genus}(C)}$. In order to generalize this to the product of
every curve of genus at least 2, we consider different
polarisations from the ones used in \cite{ross}. But in both
articles, these polarisations lie near the boundary of the ample
cone.

Later, by computing slopes, we also produce some Kodaira-fibred
surfaces with nonzero signature which is slope unstable with respect
to either the polarisations in \cite{ross} or in this paper. The
main results are the following. \vspace{0.2cm}
\\
\NI{\bf Theorem A.} Let $C$ be a smooth curve of genus $q$ at least
2. Then $X=C\times C$ is not semistable with respect to certain
polarisations. Thus $X$ admits some K\"{a}hler classes that do not
contain cscK metrics.

\vspace{0.2cm} \NI{\bf Theorem B. (Cor. \ref{nonK})} There exist
Kodaira-fibred surfaces $X$ with nonzero signature, which are not
semistable with respect to certain polarisations. Thus these
Kodaira-fibred surfaces $X$ admit some K\"{a}hler classes that do
not contain cscK metrics. \vspace{.2cm}

As consequences of Theorem A and B, we also have

\NI{\bf Corollary C.} Let $C$ bs a smooth curve of genus $q$ of
genus at least 2. Then $X=C\times C$ is not asymptotically Hilbert
semistable (resp. not asymptotically Chow semistable) with respect
to certain polarisations.

\vspace{0.2cm} \NI{\bf Corollary D.} There exist Kodaira-fibred
surfaces $X$ with nonzero signature, which are not asymptotically
Hilbert semistable (resp. not asymptotically Chow semistable) with
respect to certain polarisations.

\vspace{0.2cm} \NI{\bf Acknowledgments} I would like to thank Claude
LeBrun for suggesting this problem to me and his great guidance. I
also want to thank Masashi Ishida for reading this article and
catching many errors.

\section{slope stability}
The concept of slope stability for a polarised manifold was
introduced by Thomas and Ross in \cite{rt}, and we will put it here
for the sake of completeness of the paper. Let $(X,L)$ be a
polarised manifold, and $Z$ be a subscheme of  $(X,L)$.  If we blow
up $X$ along $Z$, we would get the exceptional divisor $E$ and the
projection map $\pi: \widehat X\to X$, which is an isomorphism
outside $E$. Now we can define the Seshadri constant of $Z$ as
$$\epsilon(Z,L) = \text{sup}\{ c\in \mathbb Q: \pi^* L-cE \text{ is ample}\}.$$
Let $n$ be the complex dimension of the manifold $X$, $x\in\QQ$ with
$kx\in \NN$, we have the Hilbert polynomials of $L$ and $\pi^*L-xE$
as $\chi(kL) = a_0k^n + a_1k^{n-1} + O(k^{n-2})$ and
$\chi(k(\pi^*L-xE)) = a_0(x)k^n + a_1(x) k^{n-1} + O(k^{n-2})$. Note
that $a_i(x)$ is a polynomial of degree at most $i$, so it can be
extended to all real $x$. Let $\tilde{a}_i(x)= a_i -a_i(x)$. For
$0<c\le \epsilon(Z,L)$, the slope of $X$ and the quotient slope of
$Z$  are defined to be
\begin{eqnarray*}
  \mu(X,L)&=&\frac{a_1}{a_0},\\
  \mu_c(\O_Z,L)&=&\frac{\int_0^c \tilde{a}_1(x) + \frac{\tilde{a}'_0(x)}{2}dx }{\int_0^c \tilde{a}_0(x)
  dx}.
\end{eqnarray*}
These are well defined since by Riemann-Roch theorem, we have
$$a_0=\frac{1}{n!}\int_X c_1(L)^n>0 ,\quad \text{and}$$
$$a_0'(x)=-\frac{1}{(n-1)!}\int_{\widehat X} c_1(L-xE)^{n-1}\cdot E<0. $$
The second equality implies $\int_0^c \tilde{a}_0(x)dx= ca_0-\int_0^c
a_0(x)dx >0 $.

\begin{defn}
  A polarised manifold $(X,L)$ is called slope semistable with respect to a
    subscheme $Z$ if
    $$\mu(X,L)\le \mu_c(\O_Z,L) \quad \text{ for all } 0<c\le
  \epsilon(Z,L).$$
 Furthermore, if it is slope semistable with respect to all subschemes, then it is called slope semistable .
\end{defn}

Since $\epsilon(Z, L^{\otimes m})=m\cdot\epsilon(Z, L)$,
$\mu(X,L^{\otimes m})=\frac{1}{m} \mu(X,L)$, and
$\mu_c(\O_Z,L^{\otimes m})=\frac{1}{m}\mu_{\frac{c}{m}}(\O_Z,L)$,
slope semistability is preserved if $L$ is replacing by its power,
and the notion of slope semistability can be extended to ample
$\QQ$-divisors.

Slope stability in fact gives an obstruction to the existence of
cscK metrics. Donaldson \cite{d}, Chen-Tian \cite{ch} and Mabuchi
\cite{mab} \cite{mab2} have made substantial progress on relating
the existence and uniqueness of extremal K\"ahler metrics in Hodge
K\"ahler classes to the K-stability of polarized projective
varieties. In particular, it has been shown that K-stability is a
necessary condition for the existence of cscK metrics for a
polarized projective variety. Ross and Thomas \cite{rt} show that
K-semistability implies slope semistability.

\begin{thm}(\cite{rt2})\label{slopeK}
  If $(X,L)$ is K-semistable, then it is slope semistable with respect to every subscheme $Z$.
\end{thm}
\begin{proof}
  See \cite{rt2} Thm.4.18.  When
  $X$ and $Z$ are smooth, this is also proved in
  \cite{rt} Thm. 4.2.
\end{proof}

Therefore, slope semistability provides an obstruction to the
existence of cscK metrics. Furthermore, since asymptotic Hilbert
(resp. Chow) semistability implies K-semistability, slope
semistability also provides an obstruction to the notions of Hilbert
(resp. Chow) semistability of a projective variety.

\subsection*{Slope semistability for smooth complex surfaces}
In this article, we only consider the case when $X$ is a smooth compact
complex surface and $Z$ is a smooth complex curve in $X$.  Since the
complex dimensions of $X$ and $Z$ are 2 and 1, the blow-up of $X$
along $Z$ would just be $X$ itself, and the Seshadri constant is
$$\epsilon(Z,L) = \sup\{ c\in \mathbb Q: L-cZ \text{ is ample} \}.$$
Let $K$ be the canonical divisor of $X$. We can express both the
slope of $(X, L)$ and the quotient slope of $Z$ in terms of the
intersection numbers (\cite{rt} Cor. 5.3):
\begin{eqnarray}\label{eq:slope}
  \mu(X,L) &=& -\frac{K\cdot L}{L^2},\\
  \mu_c(\O_Z,L) &=& \frac{3(2L\cdot Z-c(K\cdot Z+Z^2))}{2c(3L\cdot Z-cZ^2)}. \nonumber
\end{eqnarray}
Therefore the Seshadri constant $\epsilon(Z,L)$, the slope
$\mu(X,L)$, and the quotient slope $\mu_c(\O_Z,L)$ depend only on
the classes of $L$ and $Z$ modulo numerical equivalence.  And we
could extend the equations \eqref{eq:slope} to any $\QQ$-divisor $L$
even if it is not ample. Nevertheless, this might cause a zero
denominator in the computation of slopes.

\section{bounds of the ample cone}
Let $C$ be a smooth curve of genus $q$, which is grater than 1, and
$X=C\times C$. Let $p_i$ be the projection onto the $i$-th factor,
$c$ be a fixed point in $C$, and  $\gamma_i$ be the class of the
fibre $p_i^{-1}(c)$ in the N\'eron-Severi group $N^1(X)_\mathbb Q$.
The class of the canonical divisor of $X$ is
$K_X=(2q-2)(\gamma_1+\gamma_2)$ which is ample. Let $f=\gamma_1+\gamma_2$,
$\delta$ be the
class of the diagonal. For convenience, we make the change of
variables $\delta'=\delta-f$. Then we have the following
intersection numbers on $X$:
$$ f^2 =2, \quad \delta'\cdot f=0,\,\,\,\, \text{and} \quad \delta'^2=-2q.$$

\NI In this section, we consider the intersection of the ample cone
and the $f, \delta'$ plane in the N\'eron-Severi group $N^1(X)_{\QQ}$.

First of all, consider the $\QQ$-divisor
$$l_s= sf+\delta'$$
which is ample for $s\gg 0$. To find the infimum of $s$ to make
$l_s$ ample, we need the following tool: \\
\NI{\bf Nakai's criterion:} Let $X$ be an algebraic compact complex
surface, and D be a divisor on $X$. Then $D$ is ample if and only if
$D\cdot D>0$ and $D\cdot C'>0$ for each irreducible curve $C'$.

\begin{thm}\label{ls}
$l_s$ is ample if and only if $s> q$.
\end{thm}
\begin{proof}
The essential observation in the proof is the existence of  an
irreducible curve in the class $\delta=\delta'+f$. That is  the
diagonal curve $D=\{(x,x):x\in C\}\subseteq C\times C$.  If $l_s$ is
ample, we have $l_s\cdot D=2s-2q$ is positive. Therefore $s> q$, as
required. Now, suppose that $s> q$. One has $l_s\cdot l_s =
2(s^2-q)>0$, and $l_s\cdot D=2s-2q> 0$. Let $C'$ be any irreducible
curve distinct from $D$ in $X$. Since the intersection pairing is
nondegenerate, we could write
$$[C']=x_1\gamma_1+x_2\gamma_2+y\delta'+\alpha,$$
where $\alpha\in N^1(X)_\QQ$ is orthogonal to $\gamma_1, \gamma_2,$
and $\delta'$. By intersecting with $\gamma_1, \gamma_2$, we find
that $x_1, x_2 > 0$. Moreover, since $D$ and $C'$ are two
distinct irreducible curves, $[D]\cdot [C']\geq 0$, which yields
$$x_1+x_2-2qy\geq 0.$$
The direct computation shows that $l_s\cdot
[C']=s(x_1+x_2)-2gy>x_1+x_2-2gy\geq 0.$ Therefore a direct
application of Nakai's criterion implies that $l_s$ is ample.
\end{proof}
In Ross' paper \cite{ross},  he considers the case that $C$ is a
simple branched cover of $\PP_1$ of degree $k$, where $2\leq
k-1<\sqrt q$, and the $\QQ$-divisor $L_t= tf-\delta'$. Let $s_C =
\text{inf}\{ t: L_t \text{ is ample} \}.$ Kouvidakis \cite{Ku} shows
that $s_C=\frac{q}{k-1}$. Furthermore, Kouvidakis \cite{Ku} also
shows that when $C$ is a curve of general moduli, we have $\sqrt
q\leq s_C\leq \frac{q}{[\sqrt q]}$. In particular, if $q$ is a
perfect square, $s_C=\sqrt q$. By the previous discussion, we know
that every class between the thick lines in figure 1 is ample.
\begin{figure}[htbp]
  \begin{center}
    \mbox{
      \subfigure[$C$ admits a simple branched cover over $\PP_1$]{\scalebox{0.85}{\input{nef.pstex_t}}}
      \quad\quad\quad
      \subfigure[$C$ is a curve of general moduli, and $q$ is a perfect square.]{\scalebox{0.85}{\input{nef2.pstex_t}}} \quad
      }
  \end{center}
\begin{center}
Figure 1. Ample cone in $f,\delta'$ plane
\end{center}
\end{figure}

\section{unstable products of curves}
In this section, we will show that for every smooth curve of genus
$q$ greater than 1, the product $X=C\times C$ is not slope
semistable with respect to certain polarisations. Ross \cite{ross}
has shown that when $C$ admits a simple branched cover of $\PP_1$ of
degree $2\leq k-1<\sqrt q$, $X=C\times C$ is not slope semistable
with respect to the polarisations $L_t=tf-\delta'$ for $t$
sufficiently close to $s_C$. But here we consider general curves,
and the polarisations $l_s=sf+\delta'$.

\begin{thm}\label{unstableproduct}
Let $C$ be a smooth curve of genus $q$ at least 2. Then $X=C\times
C$ is not semistable with respect to the polarisations $l_s$ for $s$
sufficiently close to $q$.
\end{thm}
\begin{proof}
  By Thm. \ref{ls}, let $s>q$ so that $l_s$ is ample.  The canonical divisor of $X$ is
  $K_X=(2q-2)f$, and
  \begin{equation}
    \mu(X,l_s) = -\frac{K\cdot l_s}{l_s^2} = -\frac{s(2q-2)}{s^2-q}.\label{slopeofX}
  \end{equation}
Let $D$ be the diagonal curve, whose class is
$\delta=f+\delta'$. We now consider the Seshadri constant of $D$.
  we have that $l_s-cD = sf+\delta' -c (f+\delta') = (s-c)f+(1-c)\delta'$. And by previous discussion, it is
  ample if and only if $c< \frac{s+s_c}{1+s_c}$.  Thus $\epsilon(D,l_s)=
  \frac{s+s_c}{1+s_c}>1$. (It is obvious to see $\epsilon(D,l_s)\geq 1$ since $l_s-D=(s-1)f>0$.)
To calculate the slope of $Z$ we need the quantities:
\begin{eqnarray}
  l_s\cdot D &=& (sf+\delta').(f+\delta') \\
&=& 2s -2q, \nonumber \\
K\cdot D&=& (2q-2)f\cdot(f+\delta')= 2(2q-2),  \nonumber \\
D^2 &=&  (f+\delta')^2 =  2-2q. \nonumber
\end{eqnarray}
Thus from \eqref{eq:slope},
\begin{eqnarray}
  \mu_c(\O_D,l_s) &= &\frac{3(2l_s\cdot D-c(K\cdot D+D^2))}{2c(3l_s\cdot D-cD^2)}\nonumber  \\
  &=& \frac{3( 4s - 4q - c(2q-2))}{2c(6s-6q-2c+2cq)}. \label{slopeofD}
\end{eqnarray}
We claim that if $0<c<\frac{3}{4}$, then
$\mu_c(\O_Z,l_s)<\mu(X,l_s)$ as $s$ tends to $q$ from above.  Since
this is an open condition it is sufficient to show that it holds
when $s= q$.  By (\ref{slopeofX}, \ref{slopeofD}),
\begin{eqnarray}
\mu(X,l_q) &=& -\frac{q(2q-2)}{q^2-q}=-2,  \nonumber\\
\mu_c(\O_Z,l_q)&=&\frac{-3c(2q-2)}{2c(-2c+2cq)}=\frac{-3}{2c}
\label{slopedofD}
\end{eqnarray}
  Hence as $0<c<\frac{3}{4}$, $\mu_c(\O_Z,l_s)<\mu(X,l_s)$ as $s$
  tends to $q$ from above, which proves that $(X,l_s)$ is not slope semistable.
\end{proof}

\begin{cor}\label{nonK}
Let $C$ bs a smooth curve of genus $q$ of genus at least 2. Then
$X=C\times C$ admit some K\"{a}hler classes that do not contain cscK
metrics.
\end{cor}
\begin{proof}
It follows from (\ref{unstableproduct}) and (\ref{slopeK}).
\end{proof}

\NI This completes the proof of Theorem A in the introduction.

\begin{cor}
Let $C$ bs a smooth curve of genus $q$ of genus at least 2. Then
$X=C\times C$ is not asymptotically Hilbert semistable (resp. not
asymptotically Chow semistable) with respect to certain
polarisations.
\end{cor}
\begin{proof}
It follows from (\ref{unstableproduct}) and (\ref{slopeK}).
\end{proof}

Related to the existence of cscK metric, Mabuchi introduces Mabuchi
functional for a given K\"ahler class $\Omega$ on a complex manifold
$X$. Chen-Tian \cite{ch}and Donaldson \cite{d2} show that the
existence of cscK metrics in $\Omega$ implies that the Mabuchi
functional is bounded from below. It is conjectured that the
existence of the cscK metrics is equivalent to the properness of the
Mabuchi functional. In the case that $X$ has negative first Chern
Class, Chen \cite{c} introduces the notion of $J$-flow, and show
that the convergence of $J$-flow implies lower boundedness of the
Mabuchi functional. For a polarised surface $(X,L)$ with negative
first Chern class, Weinkove \cite{w} has the following theorem about
sufficient condition of the convergence of the $J$-flow.

\begin{thm}(\cite{w})\label{jflow}
Let $(X,L)$ be a polarised surface with negative first Chern Class. Let the divisor $\alpha$ be
defined a by
$$\alpha= 2(K.L) L - (L^2)K.$$
If $\alpha$ is ample, then the $J$-flow converges and the Mabuchi
functional is proper on the class $c_1(L)$.
\end{thm}

When $C$ is a curve of genus $q$ at least 2, $X=C\times C$ and $L=l_s=sf+\delta'$ with $s>q$
it is easy to determine when $\alpha$ is ample.

\begin{lem}\label{wample}
Let $L=l_s=sf+\delta '$ with $s>q$.
Then $\alpha=2(K.L)L-(L^2)K$ is ample if and only if $s> q + \sqrt{q^2-q}$.
\end{lem}
\begin{proof}
Since  $l_s^2 = 2s^2-2q$ and $K\cdot l_s=2s(2q-2)$, we have
\begin{eqnarray*}
    \alpha &=& 4s(2q-2)(sf+\delta') - (2s^2-2q)(2q-2)f \\
&=& 2(2q-2)( (s^2+q)f+2s\delta'),
  \end{eqnarray*}
  which is ample if and only if $s^2+q>2qs$. The conclusion follow.
\end{proof}

\begin{cor}
Let $C$ be a smooth curve of genus at least 2, and $X$ be the product $X=C\times C$.
Then the Mabuchi functional is
proper on the polarised surface $(X, l_s)$ if $s>q+\sqrt{q^2-q}$.
\end{cor}
\begin{proof}
It follows from Theorem \ref{jflow} and lemma \ref{wample}.
\end{proof}

\section{Unstable Kodaira fibrations of nonzero signature}
In this section, we give a short sketch of an explicit construction
of Kodaira fibrations with nonzero signature. Then we will show
Kodaira fibrations constructed in this way are not slope semistable
with respect to certain polarisations.

\begin{defn}
A Kodaira fibration is a holomorphic submersion $\pi:X\to B$ from a
compact complex surface $X$ to a compact complex curve $B$, with
base $B$ and fibre $F=f^{-1}(z)$ both have genus $\geq 2$.
\end{defn}
Clearly, $\pi$ is locally a trivial fibre bundle in the smooth
sense. Nevertheless, the complex structures of all fibres may vary.
A surface admitting a Kodaira fibration is called a Kodaira-fibred
surface. Every Kodaira-fibred surface $X$ is algebraic since by
the adjunction formula,
$K_X\cdot F=-\chi(F)\ge 2,$ and  $c_1^2(K_X+kF)>0$ for $k$ large
enough, where $F$ is the class of the fibre. On the other hand, $X$
could not contain any rational or elliptic curves because if $C$ is
a curve in $X$ with genus less than two, then by the fact that a
holomorphic map from a curve of lower genus to a curve of higher
genus must be constant, we have $\pi(C)$ is a point. Therefore $C$
lies in a fibre, which is absurd since the fibre has the genus
greater than or equal to 2 . The Kodaira-Enrique classification
theorem henceforth tells us that $X$ is a minimal surface of general
type. Furthermore, since $X$ contains no $(-2)$ rational curves,
we have the canonical divisor $K_X>0$.

A product $B\times F$ of two complex curves of genus $\ge 2$ is
certainly Kodaira fibred, but such a product would have signature
$\tau =0$ since it admits an orientation-reversing diffeomorphism.
We can construct Kodaira fibrations of nonzero signature in the
following way:

Let $C$ be a compact complex curve of genus $\ge 2$, $G$ be a finite
group of order divisible by $r$, which acts effectively on $C$. We
can define a homomorphism by the composition $\pi_1(C)\to H_1(C,
\ZZ)\to H_1(C, \ZZ_r)$. Since $H_1(C, \ZZ_r)$ gas finite order,
there exists an unbranched finite cover $h:B\to C$ such that
$h_*(\pi_1(B))=\text{ker}[\pi_1(C)\to H_1(C, \ZZ_r)]$. Clearly, the
genus of $B$ is greater than or equal to the genus of $C$. Let
$\Sigma\subset B\times C$ be the union of the graphs of $g\circ
h:B\to C$, where $g$ runs over $G$.

\begin{lem}\label{cyclic}
The homology class of $\Sigma$ in $H_2(B\times C,\ZZ)$ is divisible by
$r$. Hence there exists a cyclic $r$-cover  $X$ of $B\times C$
branched over $\Sigma$.
\end{lem}
\begin{proof}
(see \cite{bpv}). Since $c_1(\O(\Sigma))$ is the poincar\'{e}
duality of the fundamental class of $\Sigma$, it suffices to show
that $c_1(\O(\Sigma))$ is divisible by $r$. That is, we have to show
the intersection pairing $(c_1(\O(\Sigma)), \alpha)\equiv
0(\text{mod r})$ for all $\alpha\in H^2(B\times C, \ZZ)$. Let
$p_1:B\times C\to B$ and $p_2:B\times C\to C$ be the projections.
Using K\"{u}nneth's formula, we deal it with three cases: $\alpha
\in p_1^*(H^2(B,\ZZ))$, $\alpha \in p_2^*(H^2(C,\ZZ))$, and
$\alpha=p_1^*(H^1(B,\ZZ))\cup p_2^*(H^1(C,\ZZ))$. The first two
cases follow from the fact that $r$ divides the order of $G$. The
last case can be done by applying the projection  formula
\begin{eqnarray}
(c_1(\O(\Sigma)),p_1^*\beta\cup p_2^*\gamma) &= &(p_1^*\beta,c_1(\O(\Sigma)) \cup p_2^*\gamma ) \nonumber \\
   &= &\sum_{g\in G}(\beta,   (g\circ h)^*\gamma )\equiv 0  \,\,\,(\text{mod}\, r) \nonumber.
\end{eqnarray}
The last equality is because $h_*(\pi_1(B))=\text{ker}[\pi_1(C)\to
H_1(C, \ZZ_r)]$.
\end{proof}

The explicit construction of the surface $X$ in the lemma
\ref{cyclic} is as follows. Since The homology class of $\Sigma$ in
$H_2(B\times C)$ is divisible by $r$, we can have a line bundle $\L$
on $B\times C$ such that $\O(\Sigma)=\L^{\otimes r}$ and a section
$s\in\Gamma({B\times C, \O(\Sigma)})$, which vanishes exactly on
$\Sigma$. Let $L$ be the total space of $\L$, $\phi:L\to B\times C$
be the line bundle projection, and $t$ be the tautological section
of the pull-back bundle $\phi^*\L$ over $L$. Now take $X$ to be the
zero set of the section $\phi^*s-t^r\in\Gamma(L, \phi^*\L^{\otimes
r})$. It can be seen that $X$ is a branched $r$ to 1 cover of
$B\times C$, which is branched over $\Sigma$. And since all the
singularities of $X$ lies on the singularities of the section $s$,
and $\Sigma$ is the disjoint union of $|G|$ smooth curves, $X$ is a
smooth surface. Accordingly, $X$ inherits a natural projection
$\pi:X\to B$. This projection $\pi$ is a submersion since it admits
a local section everywhere. The signature of $X$ can be computed as
follows: Let $p=$ genus of $B$, $q=$ genus of $C$, and $d=$ degree
of $h$. By the Riemann-Hurwitz formula, we have the Euler number
\begin{eqnarray*}
\chi(X)&=&\chi(B)(r\chi(C)-(r-1)|G|) \nonumber \\
&=&4r(p-1)(q-1)+2(p-1)(r-1)|G| >0 \nonumber  \\
\end{eqnarray*}
On the other hand, we have the canonical divisor
$\it{K}_X=\pi^*\it{K}_{B \times C}+ R$, where $R$ is the
ramification divisor. To compute the self intersection number of the
canonical line bundle $K_X$, we need the following lemma.
\begin{lem}\label{ino}
Let $M$, $N$ be two compact complex manifolds of the same dimension
$m$, and $f:M\to N$ a smooth map of degree $d$. Then the self
intersection number of the graph $\Gamma$ of $f$ in $M\times N$ is given
by $d\cdot \chi(N)$.
\end{lem}
\begin{proof}
Let $p_1, p_2$ be the projection of $M\times N$ to $M$ and $N$,
respectively. Then the normal bundle $N_\Gamma$ of $\Gamma$ is isomorphic to the
pullback bundle $p_2^*T_N|_\Gamma$ of the tangent bundle of $N$
restricted to $\Gamma$ since we have the following commutative diagram between two
exact sequences:

$$\begin{CD}
 0@>>>  T_M|_\Gamma @>>> (T_M\oplus T_N)|_\Gamma  @>>>T_N|_\Gamma@>>> 0 \\
 && @VV{\thickapprox}V         @VV{\thickapprox}V    \\
 0@>>>  T_\Gamma @>>> (T_M\oplus T_N)|_\Gamma  @>>>N_\Gamma@>>> 0
\end{CD}$$
\\
Let $eu(V)$ denote the Euler class of the vector bundle $V$. We have
\begin{eqnarray*}
\Gamma\cdot \Gamma&=&\int_\Gamma eu(N_\Gamma)\\
&=& \int_\Gamma p_2^*eu(T_N) \quad\text{by naturality of Euler class}\\
&=& \int_M f^*eu(T_N)\\
&=& \int_M d\cdot\frac{\chi(N)}{\chi(M)}eu(T_M) \\
&=& d\cdot \chi(N).
\end{eqnarray*}
\end{proof}
\NI Using $K_X=\pi^*K_{B\times C}+R$, where $R=\frac{r-1}{r}\pi^*\Sigma$ is the ramification divisor, we have
\begin{eqnarray*}
\it{K_X}^\text{2}&=&rK_{B\times C}^2+2(r-1)K_{B\times C}\cdot \Sigma+\frac{(r-1)^2}{r}\Sigma\cdot \Sigma\\
&=&2r(2q-2)(2p-2)+2(r-1)((2p-2)|G|+(2q-2)d|G|)+\frac{(r-1)^2}{r}d(2-2q)|G| \nonumber \\
 &=&8r(p-1)(q-1)+4(r-1)(p-1)|G|+2\frac{r^2-1}{r}(q-1)d|G|>0 \nonumber.\\
\end{eqnarray*}
Using Signature formula, we have
$\tau(X)=\frac{1}{3}(K_X^2-2\chi(X))=
\frac{2(r^2-1)}{3r}(q-1)d|G|>0$ as a result.

\begin{eg*}
Let $C$ be a curve of genus $3$ with a holomorphic
 involution $\iota : C\to C$ without fixed points;
 one may visualize such an involution as a $180^{\circ}$ rotation
 of a $5$-holed doughnut about an axis which passes though
 the middle hole, without meeting the doughnut.
 Let $h: B\to C$ be the unique $64$-fold unbranched cover
 with $f_{*}[\pi_{1}(B)]= \ker [\pi_{1}(C) \to H_{1}(C, \ZZ_{2})]$;
 thus $B$ is a complex curve of genus $129$. Let
 $\Sigma \subset B\times C$ be the union of
 the graphs of $f$ and $\iota\circ h$. Then
 the homology class of $\Sigma$ is divisible by $2$. We may
 therefore construct a ramified double cover $X\to B\times C$
 branched over $\Sigma$. The projection $X\to B$
 is then a Kodaira fibration, with fiber $F$ of genus $6$.
 The projection $X\to C$ is also a Kodaira fibration,
 with fiber of genus $321$. The signature of this doubly Kodaira-fibres complex surface is
 $\tau (M) = 256 > 0$.
\end{eg*}

Recall that in  Ross' construction \cite{ross} of unstable products
of curves, $C$ is a compact complex curve of genus $q\geq 2$, which
admits a simple branched covering map to $\PP_1$ of degree $k$, for
$2\leq k-1<\sqrt{q}$. Let $X_0=C\times C$. Consider the
$\QQ$-divisor $\it{L}_t=tf-\delta'$, where $f=\gamma_1+\gamma_2$ is
the class of the sum of fibres in two directions, and
$\delta'=\delta -f$ with $\delta$ the class of diagonal. It is shown
in \cite{Ku} that $\it{L}_t$ is ample if $t>\frac{q}{k-1}$. Let
$Z=C\times_{\PP_1} C-\delta$ be the residual divisor of the diagonal
in the fibre product. By computing the slope of $X_0$ and $Z_0$,
Ross \cite{ross} shows that $(X_0,\it{L}_t)$ is destablized by $Z$
for $t$ sufficiently close to $\frac{q}{k-1}$ from above.

Now let $X_1$ be an unbranched $d$-covering of $X_0$
with the covering map $\pi_1:X_1\to X_0$. Let $L_{1,t}=\pi_1^* L_t$,
and $Z_1=\pi_1^*Z$, then the Seshadri constant $\epsilon(Z_1, X_1,
L_{1,t})=\epsilon(Z, X_0, L_t)$. Since $K_{X_1}=\pi_1^* K_{X_0}$,
and the intersection pairings $\pi_1^*D\cdot \pi_1^*D'=d(D\cdot D')$
for $D, D'$ any divisors on $X_0$, we have the slopes $\mu(X_1,
L_{1,t})= \mu(X_0, L_t)$, and the quotient slopes $
\mu_c(\O_{Z_1},L_{1,t})= \mu_c(\O_Z,L_t)$. It follows that
$(X_1,L_{1,t})$ is destablized by $Z_1$ for $t$ sufficiently close
to $\frac{q}{k-1}$.

Here we take $X_1=B\times C$, a $d$ to 1 unbranched cover of
$C\times C$, where $B$ and $C$ satisfy all hypotheses in the
previous construction of Kodaira fibrations. Let $B_0=B\times c$,
$C_0=b\times C$ be the classes of the fibres of the projection onto
$C$-factor and $B$- factor, respectively, which is independent of
the choices of the points $b\in B$, $c\in C$. Let $X_2$ be the
constructed Kodaira fibration, which is a cyclic cover of $X_1$ with
nonzero signature. Then we have $\pi_2:X_2\to X_1$  the covering map
branched over $\Sigma$ of degree $r$, and $\pi_1:X_1\to X_0$ the
unbranched covering map of degree $d$. Let
$f_2=(\pi_1\circ\pi_2)^*f=\pi_2^*(B_0+dC_0)$,
$\delta_2=(\pi_1\circ\pi_2)^*\delta=\pi_2^*(\text{graph of}\,\, h)$.
We see that $f_2^2=2rd$, $f_2\cdot \delta_2=2rd$,
$\delta_2^2=rd(2-2q)$. For convenience make the change of variables
$\delta'_2=\delta_2-f_2$. Then we have the following intersection
numbers on $X_2$:
$$ f_2^2 =2rd, \quad \delta_2'.f_2=0,\,\,\,\, \text{and} \quad \delta_2'^2=-2rdq.$$

\NI Now consider the $\mathbb Q$-divisor
$$ L_{2,t,\varepsilon} = (\pi_1\circ\pi_2)^*L_t+\varepsilon K_{X_2}= tf_2-\delta_2'+\varepsilon K_{X_2}$$
which is ample for $\varepsilon >0$, and $t\gg 0$. Here the
canonical divisor is $K_{X_2}=\pi_2^*(K_{X_1})+R$ with R the
ramification divisor.   We define
$$ s_{\varepsilon} = \text{inf}\{ t: L_{2,t,\varepsilon} \text{ is ample} \}.$$
Clearly $s_{\varepsilon}\leq\frac{q}{k-1}$ for $L_{1,\frac{q}{k-1}}$ is
numerically effective (see \cite{Ku}) and $K_{X_2}$ is ample.

\begin{thm}\label{unstable kodaira}
Let the divisor $Z_2$ be defined by
$Z_2=\pi_2^*Z_1=(k-1)f_2-\delta_2'$. Then $X_2$ is not slope
semistable with respect to $Z_2$ for the polarisations
$L_{2,t,\varepsilon}$ if $t$ is sufficiently close to
$\frac{q}{k-1}$, and $\varepsilon$ is small enough.
\end{thm}
\begin{proof}
Let $t> s_\varepsilon$ so $L_{2,t,\varepsilon}$ is ample for any
$\varepsilon>0$. The canonical divisor of $X_2$ is
$K_{X_2}=\pi_2^*K_{X_1}+R$, where $R$ is the ramification divisor,
and $\pi_2^*\Sigma=\frac{r}{r-1}R$. From Lemma \ref{ino}, we can
compute the intersection numbers
\begin{eqnarray}
R\cdot\pi_2^*L_{1,t}&=&r(\frac{r-1}{r}\Sigma)\cdot L_{1,t}\nonumber\\
&=&(r-1)\Sigma\cdot ((t+1)(dC_0+B_0)-\text{graph of}\,\,h)\nonumber\\
&=&2(r-1)d((t+1)|G|+q-1),\nonumber\\
R\cdot Z_2&=& (r-1)\Sigma\cdot ((k-1)f_2-\delta'_2)\nonumber\\
&=& (r-1)(2kd|G|-d(2-2q))\nonumber\\
&=&2d(r-1)(k|G|-1+q), \label{rz}
\end{eqnarray} and we have
\begin{eqnarray}
\mu(X_2,L_{2,t,\varepsilon})& = & -\frac{K_{X_2}\cdot
L_{2,t,\varepsilon}}{L_{2,t,\varepsilon}^2} =
-\frac{(\pi_2^*K_{X_1}+R)\cdot(\pi_2^*L_{1,t}+\varepsilon K_{X_2})}{
(\pi_2^*L_{1,t}+\varepsilon
K_{X_2})^2} \nonumber\\
&=&-\frac{rK_{X_1}\cdot L_{1,t}+ R\cdot \pi_2^*L_{1,t}+O(\varepsilon)
}{rL_{1,t}^2+O(\varepsilon)}\nonumber \\
&=&-\frac{K_{X_1}\cdot L_{1,t}}{L_{1,t}^2}-\frac{2(r-1)d((t+1)|G|+q-1)}{2rd(t^2-q)}+O(\varepsilon)\nonumber \\
&=& \mu(X_1,
L_{1,t})-\frac{(r-1)((t+1)|G|+q-1)}{r(t^2-q)}+O(\varepsilon)\label{eq:slopeofX_2}.
\end{eqnarray}
Recall that
$L_t-Z=(t-k+1)f$ is ample when $2\leq k-1<\sqrt{q}$ (see \cite{ross}).
Therefore $L_{2,t,\varepsilon}-Z_2=(\pi_1\circ
\pi_2)^*(L_t-Z)+\varepsilon K_{X_2}$ is ample , and $\epsilon(Z_2,
L_{2,t,\varepsilon})\geq 1$ for any positive $\varepsilon$.

To calculate the quotient slope of $Z_2$, we have from
\eqref{eq:slope} and (\ref{rz}),
\begin{eqnarray}
  \mu_1(\O_{Z_2},L_{2,t,\varepsilon}) &= &\frac{3(2L_{2,t,\varepsilon}\cdot Z_2-(K_{X_2}\cdot Z_2+Z_2^2))}
  {2(3L_{2,t,\varepsilon}\cdot Z_2-Z_2^2)} \nonumber \\
  &=& \frac{3(2rL_{1,t}\cdot Z_1-rK_{X_1}\cdot Z_1-rZ_1^2)-3R\cdot \pi_2^*Z_1+O(\varepsilon)}{2(3rL_{1,t}\cdot Z_1-rZ_1^2)+O(\varepsilon)}
  \nonumber.\\
  &=& \mu_1(\O_{Z_1},L_{1,t})-\frac{3R\cdot \pi_2^*Z_1}{2r(3L_{1,t}\cdot Z_1-Z_1^2)}+O(\varepsilon) \nonumber.\\
 &=& \mu_1(\O_{Z_1},L_{1,t})-\frac{6d(r-1)(k|G|+q-1)}{4rd(3t(k-1)-(k-1)^2-2q)}+O(\varepsilon) \nonumber.\\
&=&
\mu_1(\O_{Z_1},L_{1,t})-\frac{3(r-1)(k|G|+q-1)}{2r(3t(k-1)-(k-1)^2-2q)}+O(\varepsilon)
\label{eq:slopeofZ_2}.
\end{eqnarray}
Since $\mu_1(\O_{Z_1},L_{1,t})<\mu(X_1,L_{1,t})$ near
$t=\frac{q}{k-1}$ from Ross' computation \cite{ross}, by
(\ref{eq:slopeofX_2}) and (\ref{eq:slopeofZ_2}), it suffices to show
$$\frac{(r-1)((t+1)|G|+q-1)}{r(t^2-q)}<
\frac{3(r-1)(k|G|+q-1)}{2r(3t(k-1)-(k-1)^2-2q)}. $$  Because the
slopes depend continuously on $t$,  we just need to show the
inequality holds at $t=\frac{q}{k-1}$. That is,
$2\frac{(k-1)^2}{q}((\frac{q}{k-1}+1)|G|+q-1)<3(k|G|+q-1)$. Using
the assumption $2\leq k-1<\sqrt q$, we have
\begin{eqnarray*}
&&3(k|G|+q-1)-2\frac{(k-1)^2}{q}((\frac{q}{k-1}+1)|G|+q-1) \\
&=&3(k|G|+q-1)-2((k-1+\frac{(k-1)^2}{q})|G|+\frac{(k-1)^2(q-1)}{q})\\
&>&3(k|G|+q-1)-2((k-1+1)|G|+(q-1))\\
&=&k|G|+q-1>0.
\end{eqnarray*}
Thus $\epsilon(Z_2,L_{2,t,\varepsilon})\ge 1$ and
$\mu_1(\O_{Z_2},L_{2,t,\varepsilon})<\mu(X_2,L_{2, t,\varepsilon})$
as  $t$ is sufficiently close
to $\frac{q}{k-1}$, and $\varepsilon$ tends to $0$ from above, which proves that $(X_2,L_{2,t,\varepsilon})$ is
not slope semistable.
\end{proof}

Instead of the polarisations
$L_{2,t,\varepsilon}=tf_2-\delta'_2+\varepsilon K_{X_2}$, we could
also consider the polarisations
${l_2,s,\varepsilon}=sf_2+\delta'_2+\varepsilon K_{X_2}$, where
$\varepsilon>0$. Let $D$ be the diagonal in $X=C\times C$,
$D_1=\pi_1^*D$, and  $D_2=\pi_2^*D_1$.

\begin{thm}\label{unstable kodaira2}
$X_2$ is not slope semistable with respect to the curve $D_2$ for
the polarisations $l_{2,s,\varepsilon}$ if $s$ is sufficiently close
to $q$, and $\varepsilon$ is small enough.
\end{thm}
\begin{proof}
Let $s> q$ so that $l_{2,s,\varepsilon}$ is ample for any
$\varepsilon>0$. The canonical divisor of $X_2$ is
$K_{X_2}=\pi_2^*K_{X_1}+R$, where $R$ is the ramification divisor.
From (\ref{ino})and (\ref{slopedofD}), we have
\begin{eqnarray*}
R\cdot\pi_2^*l_{1,s}&=&r(\frac{r-1}{r}\Sigma)\cdot l_{1,s}\nonumber\\
&=&(r-1)\Sigma\cdot ((s-1)(dC_0+B_0)+\text{graph of }\,h)\nonumber\\
&=&2(r-1)d((s-1)|G|-q+1),\nonumber\\
\end{eqnarray*}
and
\begin{eqnarray}
\mu(X_2,l_{2,s,\varepsilon})& = & -\frac{K_{X_2}\cdot
l_{2,s,\varepsilon}}{l_{2,s,\varepsilon}^2} =
-\frac{(\pi_2^*K_{X_1}+R)\cdot(\pi_2^*l_{1,s}+\varepsilon K_{X_2})}{
(\pi_2^*l_{1,s}+\varepsilon
K_{X_2})^2} \nonumber\\
&=&-\frac{rK_{X_1}\cdot l_{1,s}+R\cdot \pi_2^*l_{1,s}+O(\varepsilon)
}{rl_{1,s}^2+O(\varepsilon)}\nonumber \\
&=& \mu(X_1, l_{1,s})-\frac{2(r-1)d((s-1)|G|-q+1)}{2rd(s^2-q)}+O(\varepsilon)\nonumber \\
&=&
\frac{-s(2q-2)}{s^2-q}-\frac{(r-1)((s-1)|G|-q+1)}{r(s^2-q)}+O(\varepsilon)\label{slopeofX_2}.
\end{eqnarray}
To bound the Seshadri constant of $D_2$, we have $\epsilon(
l_{2,s,\varepsilon}, D_2)\geq 1$ since
$l_{2,s,\varepsilon}-D_2=(\pi_1\circ \pi_2)^*l_s+\varepsilon
K_{X_2}$ is ample if $s>q$. To calculate the quotient slope of
$D_2$, we need the quantities:
\begin{eqnarray}
  l_{2,s,\varepsilon}\cdot D_2 &=& r(sf_1+\delta_1')(f_1+\delta_1')+O(\varepsilon) \nonumber\\
&=& r(2ds-2dq)+O(\varepsilon), \nonumber \\
R\cdot D_2&=& (r-1)\Sigma(f_1+\delta_1')= (r-1)d(2-2q),  \nonumber \\
 D_2^2 &=&  rD_1^2=  rd(2-2q) . \nonumber
\end{eqnarray}
Thus from \eqref{eq:slope} and (\ref{slopedofD}),
\begin{eqnarray}
  \mu_c(\O_{D_2},l_{2,s,\varepsilon}) &= &\frac{3(2l_{2,s,\varepsilon}\cdot D_2-c(K_{X_2}\cdot D_2+D_2^2))}
  {2c(3l_{2,s,\varepsilon}\cdot D_2-cD_2^2)} \nonumber \\
  &=& \mu_c(\O_{D_1},l_{1,s})-\frac{3cR\cdot D_2}{2c(3l_{2,s,\varepsilon}\cdot D_2-cD_2^2)}+O(\varepsilon) \nonumber.\\
 &=& \frac{3(4s-4q-c(2q-2))}{2c(6s-6q-2c+2cq)}-\frac{3cd(r-1)(2-2q)}{2rc(6ds-6dq-cd(2-2q))}+O(\varepsilon)\nonumber \\
 &=&  \frac{3(4s-4q-c(2q-2))}{2c(6s-6q-2c+2cq)}- \frac{3(r-1)(1-q)}{2r(3s-3q-c(1-q))}+O(\varepsilon) \label{slopeofD_2} .
\end{eqnarray}

We claim that
$\mu_c(\O_{D_2},l_{2,s,\varepsilon})<\mu(X,l_{2,s,\varepsilon})$ as
$s$ tends to $q$ from above. Since this is an open condition, it
suffices to show it holds when $s=q$. By (\ref{slopeofX_2}) and
(\ref{slopeofD_2}),
\begin{eqnarray*}
\mu(X_2,l_{2,q,\varepsilon}) &=& -2-\frac{(r-1)(|G|-1)}{rq}+O(\varepsilon), \\
\mu_c(\O_{D_2},l_{2,q,\varepsilon})&=&\frac{-3}{2c}+\frac{3(r-1)}{2cr}+O(\varepsilon).
\end{eqnarray*}
Hence as $c<\frac{3q}{4q+2(r-1)(|G|-1)}<\frac{3}{4}$, $s$ close to $q$, and
$\varepsilon$ small enough,
$\mu_c(\O_{D_2},l_{2,s,\varepsilon})<\mu(X_2,l_{2, s,\varepsilon})$,
which proves that $(X_2,l_{2,s,\varepsilon})$ is not slope
semistable.
\end{proof}

\begin{cor}\label{nonK}
There exist Kodaira-fibred surfaces $X$ with nonzero signature,
which admit some K\"{a}hler classes that do not contain cscK
metrics.
\end{cor}
\begin{proof}
It follows from (\ref{unstable kodaira}), (\ref{unstable kodaira2}),
and (\ref{slopeK}).
\end{proof}

\begin{cor}
There exist Kodaira-fibred surfaces $X$ with nonzero signature,
which are not asymptotically Hilbert semistable (resp. not
asymptotically Chow semistable) with respect to certain
polarisations.
\end{cor}
\begin{proof}
It follows from (\ref{unstable kodaira}), (\ref{unstable kodaira2}),
and (\ref{slopeK}).
\end{proof}

\smallskip
\NI Department of Mathematics,\\
Stony Brook University,\\
Stony Brook, NY 11794, USA.
{\small  \newline \noindent\noindent {\tt yjshu@math.sunysb.edu}}
\end{document}